\input amstex
\documentstyle{amsppt}
\document
\loadmsam
\magnification 1200
\pageheight{7.25 in}
\define \x{\tilde b}
\define \X{\tilde B}

\define \modp{\mod p}
\define \zgp{\zeta_{G,p}^\triangleleft}
\define \zlp{\zeta_{L,p}^\triangleleft}
\define \aaa{\alpha_1}
\define \ab{\alpha_2}
\define \ac{\alpha_3}

\define \pn{p^N}
\define \pna{p^{N_1}}
\define \y{c\pna}
\define \pnb{p^{N_2}}
\define \pnc{p^{N_3}}

\define \pnab{p^{N_1+N_2}}
\define \pnac{p^{N_1+N_3}}
\define \pnbc{p^{N_2+N_3}}
\define \pnabc{p^{N_1+N_2+N_3}}
\define \modpab{\text{mod }\pnab}
\define \modpbc{\text{mod }\pnbc}
\define \modpabc{\text{mod }\pnabc}
\define \mg{\pmatrix \aaa X+\ab Z & X &   Y+\ac Z \\  X & Z & 0 \\ Y & 0 &
X  \endpmatrix}
\define \ma{\pmatrix \aaa & 0 & \ab \\ 1 & 0 & 0 \\  0 &  1 & \ac
\endpmatrix}
\define \mb{\pmatrix 1 & 0 & 0 \\ 0 & 0 & 1 \\ 0 & 0 & 0 \endpmatrix}
\define \mc{\pmatrix 0 & 1 & 0 \\ 0 & 0 & 0 \\ 1 & 0 & 0 \endpmatrix}
\define \md{\pmatrix \aaa & 0 & \ab \\ 1 & 0 & 0 \\ 0 & 1 & 0 \endpmatrix}
\define \me{\pmatrix 1 & 0 & 0 \\ 0 & 0 & 1 \\ 0 & 0 & 0 \endpmatrix}
\define \mf{\pmatrix 0 & 1 & \ac \\ 0 & 0 & 0 \\ 1 & 0 & 0 \endpmatrix}
\define\zp{\Bbb Z_p}

\define\bbbz{\Bbb Z}

\define\r {\Cal R}
\define\la {\lambda}

\nologo
\topmatter
\title Associating Curves of Low Genus to Infinite Nilpotent Groups via the 
zeta function \endtitle
\author Cornelius Griffin \endauthor
\affil School of Mathematical Sciences \\ University of Nottingham \\ NG7 
2RD \endaffil
\address{Department of Mathematics, The University of Nottingham, NG7 2RD 
Nottingham, England. Email: pmxcjg\@maths.nott.ac.uk} \endaddress
\date{August 30, 2002} \enddate
\abstract It is known from work of du Sautoy and Grunewald in \cite{duSG1} 
that the zeta functions counting subgroups of finite index in infinite 
nilpotent groups depend upon the behaviour of some associated system of 
algebraic varieties on reduction $\modp.$ Further to this, in \cite{duS3, 
duS4} du Sautoy constructed a group whose local zeta function was determined 
by the number of points on the elliptic curve $E:Y^2=X^3-X.$ In this work we 
generalise du Sautoy\rq s construction to define a class of groups whose 
local zeta functions are dependent upon the number of points on the 
reduction of a given elliptic curve with a rational point. We also construct 
a class of groups that behave the same way in relation to any curve of genus 
2 with a rational point. We end with a discussion of problems arising from 
this work.
\endabstract
\endtopmatter
\leftheadtext{}
\rightheadtext{}

\heading{ \bf{Section 1: Introduction}}
\endheading
\vskip 0.25in

In \cite{GSS} Grunewald, Segal and Smith introduced the notion of a zeta 
function for an infinite group $G$ encoding (normal) subgroups of finite 
index:
$$
\zeta_{G}(s) =\sum_{H<_f G} |G:H|^{-s}
$$
and
$$
\zeta_{G}^{\triangleleft} (s) =\sum_{H\triangleleft_f G} |G:H|^{-s}.
$$
In particular they considered these functions for an infinite nilpotent 
group, as for groups of this type, the global zeta functions split as an 
Euler product of local zeta functions:
$$
\zeta_G^{(\triangleleft)}(s) = \prod_{p \text { 
prime}}\zeta_{G,p}^{(\triangleleft)}(s)=\prod_p \sum_{H<_pG}|G:H|^{-s}.
$$
Since then most effort in this subject area has gone into understanding the 
nature of these local factors for specific torsion free nilpotent groups. In 
\cite{GSS} it was shown that there exists a Lie algebra $L$ over $\bbbz$ 
associated to $G$ so that for almost all primes we have
$$
\zeta_{L,p}^{(\triangleleft)}(s) =\zeta_{G,p}^{(\triangleleft)}(s)
$$
where $\zeta_{L,p}^{(\triangleleft)}(s)$ counts subalgebras (ideals) of 
finite index in $L.$

Furthermore the authors demonstrated that  
$\zeta_{L,p}^{(\triangleleft)}(s)$ could be expressed as a $p$-adic integral 
over $Tr_d(\zp)$, the upper triangular $d\times d$ matrices over $\zp$, and 
by applying some model theory established the rationality in $p^{-s}$ of 
these functions.

By evaluating the integrals explicitly, du Sautoy and Grunewald demonstrated 
an intriguing link between the zeta function and the arithmetic of some 
algebraic varieties. In particular they showed that the zeta function of an 
infinite nilpotent torsion free group is dependent upon the number of points 
on the reduction $\mod p$ of some associated system of  algebraic varieties. 
The question then is: what type of varieties can arise in the evaluation of 
the zeta function of an infinite torsion free nilpotent group?

Du Sautoy provided the first interesting answer \cite{duS3, duS4} to this 
question by constructing a group $G(E)$ for which we have
$${}\zgp (s) =P_1(p,p^{-s})+|E(\Bbb F_p)|P_2(p,p^{-s})
$$
for rational functions $P_1,P_2$ and for the elliptic curve $E:Y^2=X^3-X.$ 
This provided the first example of  a nilpotent group with a non-finitely 
uniform zeta function.

The aim of this paper is to extend this work of du Sautoy and produce a 
larger class of algebraic varieties whose reduction $\mod p$ is encoded in 
the subgroup structure of some infinite nilpotent group. In particular we 
prove
\proclaim {Theorem 1}
Let $E$ be an elliptic curve defined over $\Bbb Q$ with a rational point in
$\Bbb Q.$ Then there exists a 9 generated, class 2 infinite torsion-free
nilpotent group $G$, associated Lie algebra $L,$ associated lines $M_1,M_2$ 
and
rational functions $P_1,\dots ,P_5 \in \Bbb Q(X,Y)$ so that for almost all 
primes
$p$, in particular  including primes dividing neither the discriminant nor 
the
coefficients of the curve, we have
$$
\split
\zgp (s) = \zlp (s) & = P_1(p,p^{-s}) +|E(\Bbb F_p)| P_2(p,p^{-s})+|M_1\cap
E(\Bbb F_p)|P_3(p,p^{-s}) \\ & +P_4(p,p^{-s})|M_2 \cap E(\Bbb F_p)| 
+P_5(p,p^{-s})|M_1\cap M_2 \cap E(\Bbb F_p)|.
\endsplit
$$
Furthermore $P_2 \not\equiv 0.$
\endproclaim
\proclaim{Theorem 2}
Let $C$ be a curve of genus 2 over $\Bbb Q$ with a rational point in $\Bbb
Q.$ Then there exists a 15 generated, class 2 torsion-free nilpotent group G
and associated Lie algebra L, $l \in \Bbb N,$ rational functions $P, Q_1,
\dots, Q_l \in \Bbb Q(X,Y)$ and varieties $V_1, \dots, V_l$ defined over 
$\Bbb Q$ so that for
almost all primes $p$
$$
\zeta_{G,p}^{\triangleleft}(s) =\zeta_{L,p}^{\triangleleft} (s) = |C(\Bbb
F_p)|P(p,p^{-s}) + \sum_{i=1}^l |V_i(\Bbb F_p)| Q_i(p,p^{-s}).
$$
Furthermore it is strictly necesary to count points on the curve $C$ in the 
evaluation of the zeta function. In particular, the $V_i$ are varieties of 
genus smaller than 2 and $C$ occurs in the subring of the motivic zeta 
function one can associate to the group $G.$

\endproclaim

The method of proof is as follows: given a torsion free nilpotent group with 
a presentation
$$
G=\langle X_1, \dots , X_d:[X_i,X_j] = \prod_{k=1}^dX_k^{a_{ij}^k} \rangle .
$$
we take the Lie algebra $L$ associated to $G$ via the Mal\rq cev 
correspondence which has a presentation
$$
L=\langle e_1, \dots , e_d:(e_i,e_j) = \sum_{k=1}^d a_{ij}^ke_k \rangle .
$$

Defining $C_j$ for $j=1, \dots, d$ to be
the matrices with $(i ,k)$-entry $c_{ik}(j)$ where
$$
(e_i,e_j) = \sum_{k=1}^d c_{ik}(j)e_k.
$$
it is known that \cite {duSG1}
$$
\zlp (s) = (1-p^{-1})^{-d} \int_{V_p^\triangleleft} |m_{11}|^{s-1}\dotsb
|m_{dd}|^{s-d} |dx|
$$
where here we define
$$
V_p^\triangleleft = \{ M \in Tr_d(\zp) :  \underline {m_i}C_jM^+ = m_{11}
\dotsb m_{dd}(Y_{ij}^1, \dots , Y_{ij}^d) \text { for some }
Y_{ij}^k \in \zp \}
$$
and $|dx|$ is the normalized Haar measure on $Tr_d(\zp).$ Here we have
denoted by $Tr_d(\zp)$ the $d\times d$ upper triangular matrices with
entries from $\zp.$ We then evaluate this integral by parts.

The paper is organised as follows: we prove Theorem 1 in Sections 2 and 3, 
and then Theorem 2 in Section 4. In Section 5 we discuss the associated 
problems of evaluating the zeta functions attached to the groups in question 
that count all subgroups of finite index, not merely normal subgroups. In 
Section 6 we discuss some problems arising from this work. We include in an 
Appendix the determinants arising in the calculation of the zeta functions 
in Sections 2 and 3.

This work is part of the author\rq s PhD thesis, Nottingham 2002, carried 
out with the support of an EPSRC studentship, and under the supervision of 
Professor Ivan Fesenko. The author would like to thank EPSRC for their 
financila support, and also I. Fesenko, M. du Sautoy and M. Edjvet for their 
helpful comments.

\vskip 0.5in
\heading {\bf {Section 2: Proof of Theorem 1}}
\endheading
\vskip 0.25in

Let $E$ be an elliptic curve defined over $\Bbb Q$ with a rational point in
$\Bbb Q.$ By means of a linear shift in $X,Y$ we may assume that $0 \in
E(\Bbb Q)$ and so $E$ has an equation of the form \cite{M1}
$$
Y^2+a_1Y+a_2XY=X^3+a_3X^2+a_4X.
$$
A transformation $Y \mapsto Y-(a_2/2)X$ enables us to write the curve as
$$
Y^2+\ac Y=X^3+\aaa X^2+\ab X
$$
or projectively as
$$
Y^2Z+\ac YZ^2=X^3+\aaa X^2Z+\ab XZ^2.
$$
Notice that this curve may be expressed as the determinant of the following
matrix:
$$
F= (f_{ij}) := \mg
$$
for appropriate $\alpha_i \in\bbbz.$ We define the Lie algebra $L$ to be
$$
L= \langle A_1, \dots, A_3,B_1, \dots, B_3, X,Y,Z : (A_i,B_j)= f_{ij}(X,Y,Z)
\rangle
$$

So how does the calculation of the zeta function
associated
to this Lie algebra differ from that presented in
\cite{duS3}? The simple answer is: not a lot. The working is made more
difficult
due to the fact that the matrix in this case is not
symmetric and so a lot of details that could be brushed under the carpet
previously now have to be confronted head on. Also the measure of sets that
we need
to calculate to show the dependence on the curve is more
difficult to realise. In any case I will
now go on to give the calculations in full. Notice that these calculations
are only valid when we consider the local zeta function of $L$
at primes $p$ not dividing  non-zero members of the set $\{\aaa , \ab, \ac
\}$ and also not
dividing the discriminant of the elliptic curve.

As outlined above we may write the zeta function as an integral
$$
\zlp (s) = (1-p^{-1})^{-9} \int_{V_p^\triangleleft} |m_{11}|^{s-1}\dotsb
|m_{99}|^{s-9} |dx|
$$
However the algebra we are working with is class 2 and so we may rewrite
this integral as $$
\zlp (s) = (1-p^{-1})^{-9} \int_{W_p^\triangleleft} |m_{11}|^{s-1}\dotsb
|m_{66}|^{s-6} |n_1|^{s-7}\dotsb |n_3|^{s-9} |dm|\cdot |dn|
$$
where now $dm$ and $dn$ are respectively the additive Haar measures on
$Tr_6(\zp)$ and $Tr_3(\zp),$ and $W_p^\triangleleft$ consists of
pairs of matrices

$$
(M,N) \in Tr_6(\zp) \times Tr_3(\zp)
$$
so that for $j=1,2,3$ we have
$$
(m_{i4},m_{i5},m_{i6})C(j)N^+ = (\alpha_0 ,\beta _0,\gamma _0)n_1n_2n_3
$$
whereas for $j=4,5,6$
$$
(m_{i1},m_{i2},m_{i3})C(j)N^+ = (\alpha_0 ,\beta_0,\gamma_0)n_1n_2n_3.
$$
Here $\alpha_0, \beta_0, \gamma_0 \in \zp$,
$$
N = \pmatrix  n_1 & a & b \\ 0 & n_2 & c \\ 0 & 0 & n_3 \endpmatrix ,
$$
and
$$
C(1) = \ma , \quad C(2) = \mb , \quad C(3) = \mc
$$

$$
C(4) = \md , \quad C(5) = \me , \quad C(6) = \mf .
$$
In other words, we are integrating by parts: we fix a basis for the centre
of the algebra and count bases for the abelianisation lying above this
particular central basis. Then we count occurrences of bases for the centre.
We may now write the zeta function as a sum
$$
\zlp (s) = \sum \Sb M_1, \dots, M_6, \\ N_1, \dots, N_3 \in \Bbb N \endSb
p^{-M_1s}\dotsb p^{-N_3(s-8)} \mu (M_1, \dots, N_3)
$$
where now $\mu (M_1, \dots , N_3)$ is the measure of those matrices $(M,N)$
with $p^{M_i}, p^{N_i}$ replacing $m_i,n_i.$
So evaluating the sum now reduces to the problem of calculating the measure
of this set. It is in this measure that the elliptic curve and
associated lines will appear. The measure can again be written as a $p$-adic
integral
$$
\mu (M_i, N_i) = \int\limits_{(a,b,c)\in \zp} \mu(\Omega_1)  \dotsb \mu
(\Omega_6) \left( \sum_{m\geq 1} \mu(\Omega_1)(p^{-m}-p^{-m-1})\right)^3 
|da|\cdot |db| \cdot |dc|
$$
where 1) $\Omega_1$ is the set of $(m_2,m_3) \in \zp^2 $ so that for
$j=4,5,6$ there
exists $(\alpha_0, \beta_0, \gamma_0) \in \zp^3$ so that
$$
(p^{M_1},m_2,m_3)C(j)N^+ = (\alpha_0,\beta_0,\gamma_0) \pnabc,
$$

2) $\Omega_2 $ is the set of $m_3\in \zp$ so that for $j=4,5,6$ there exists
$(\alpha_0, \beta_0, \gamma_0) \in \zp^3$ so that
$$
(0,p^{M_2},m_3)C(j)N^+ = (\alpha_0,\beta_0,\gamma_0) \pnabc,
$$

and 3)
$$
\Omega_3  = \cases 1 & \text { if }(0,0,p^{M_3})C(j)N^+=
(\alpha_0,\beta_0,\gamma_0)\pnabc  \\
                             0 &  \text{ otherwise.}
\endcases
$$
We  can similarly define $\Omega_4, \dots, \Omega_6$ as follows:

4) $\Omega_4$ is the set of $(m_5,m_6) \in \zp^2 $ so that for $j=1,2,3$
there
exists $(\alpha_0, \beta_0, \gamma_0) \in \zp^3$ so that
$$
(p^{M_4},m_5,m_6)C(j)N^+ = (\alpha_0,\beta_0,\gamma_0) \pnabc,
$$

5) $\Omega_5 $ is the set of $m_6\in \zp$ so that for $j=1,2,3$ there exists
$(\alpha_0, \beta_0, \gamma_0) \in \zp^3$ so that
$$
(0,p^{M_5},m_6)C(j)N^+ = (\alpha_0,\beta_0,\gamma_0) \pnabc,
$$

and 6)
$$
\Omega_6 =\cases  1 & \text{ if }
(0,0,p^{M_6})C(j)N^+= (\alpha_0,\beta_0,\gamma_0)\pnabc  \\
                             0&  \text{ otherwise.}
                             \endcases
$$

Following the notation of du Sautoy we set $\x := ac-b\pnb$ and then
one can check that investigating the above conditions, we can rewrite them 
as follows:

\proclaim{Evaluating $\mu(\Omega_i)$}
\endproclaim
To calculate
the value of $\mu(\Omega_1)$ notice that the conditions for a point to be in 
the
set become
$$
(\aaa p^{M_1+N_2+N_3}+ m_2\pnbc ,p^{M_1+N_2+N_3},  m_3 \pnbc ) \equiv
0 \modpabc
$$
$$
(p^{M_1},m_2,m_3)\pmatrix -a\aaa & -a & \pna \\ -a & 0 & 0 \\
\pna & 0 & -a \endpmatrix \equiv   0 \modpab
$$
and
$$
(p^{M_1},m_2,m_3)\pmatrix \aaa \x+\ab \pnab & \x & - \y +\ac\pnab
\\
\x & \pnab & 0 \\ -\y& 0 &  \x \endpmatrix
$$
is $0\mod p^{N_1+N_2+N_3}.$

Thus we have
\proclaim{Lemma 2.1} $\mu (\Omega_1)$ is given by:
$$
\mu(\Omega_1) = \cases 0,&\text { if } M_1<N_1 \\
p^{-2N_1}\mu(\Omega^{\prime}), &\text { if } M_1 \geq N_1
\endcases
$$
where $\Omega_1^{\prime}$ is the set of those $(m_2,m_3)\in \zp^2$ so that
$$
(p^{M_1},m_2,m_3)(S_1,S_2) \equiv   0 \modpbc
$$
and
$$
\split
(S_1,S_2)=& \biggl( \matrix -a\aaa \pnc & -a\pnc &   \pnac  \\
-a \pnc & 0 & 0 \\
\pnac & 0 & -a \pnc
\endmatrix
\cdots \\
& \cdots
\matrix
\aaa \x +\ab\pnab & \x & - \y+\ac\pnab  \\
\x & \pnab & 0 \\
-\y  & 0 &  \x
                    \endmatrix \biggr)
                    \endsplit
$$
\endproclaim
So given a solution
$(p^{M_1},X,Y)$ to this congruence, all other solutions will be
of the form $(p^{M_1},X,Y)+(0,m_2,m_3)$ where here $(m_2,m_3)$ is a solution
to the congruence
$$
\split
(m_2,m_3)&
\biggl( \matrix -a \pnc & 0 & 0 \\
\pnac & 0 & -a \pnc \endmatrix  \cdots  \\
& \qquad \cdots \matrix
\x & \pnab & 0 \\
-\y & 0 &  \x
\endmatrix \biggr)  \equiv   0  \modpbc
\endsplit \tag 2.1
$$
So to sum up, the value for $\mu(\Omega_1)$ is contained in the following
\proclaim{Proposition 2.2} If  $M_1<U_1+U_2-W_1-W_2-W_3+N_1+N_2+N_3$ then 
$\mu(\Omega_1) =0.$ For all other values of $M_1$ we have
$$
\mu(\Omega_1) = p^{U_1+U_2-(N_1+N_2+N_3)}.
$$
\endproclaim
Here we have defined
$$
U_1 := \min \bigl\{ u_1,N_2+N_3\bigr\}
$$
and
$$
U_2 :=\min \bigl\{u_2,N_2+N_3\bigr\}
$$
with
$$
u_1:=\min \bigl \{ v(detX):X \text{ a  } \, 1\times 1\text{ minor of matrix
in
\thetag{2.1}}\bigr \}
$$
and
$$
u_2:=\min \bigl \{v(detX):X \text{ a  } \, 2\times 2\text{ minor of  matrix 
in \thetag{2.1}} \bigr \}-u_1.
$$
We can evaluate these by hand and so find that
$$
u_1=\min  \bigl \{ v(a)+N_3, N_1+N_3, v(\x),N_1+N_2, v(c)+N_1 \bigr \}
$$
and
$$
\split
u_2= & \min \bigl\{ 2v(a)+2N_3, v(b)+N_1+N_2+N_3, 
2N_1+N_2+N_3,v(a)+v(\x)+N_3, \\
& v(a)+N_1+N_2+N_3, v(c)+2N_1+N_2, 2v(\x), v(\x)+N_1+N_2  \bigr\}-u_1.
\endsplit
$$
Similarly we have set
$$
W_i =\min \bigl\{N_2+N_3,w_i\bigr\}
$$
with the $w_i$ defined as the $u_i$ were above. Thus one may check, using 
the determinants evaluated in the Appendix to this paper that we have the 
following values for $w_1,w_2,w_3:$
$$
\split
w_1= & \min \bigl \{  v(a)+N_3, N_1+N_3, v(\x), N_1+N_2, v(c)+N_1, \\
& v(\ac\pnb -c)+N_1, v( \aaa\x+\ab \pnab)
\bigr \};
\endsplit
$$
$$
\split
w_2 =  \min \bigl \{ & 2v(a)+2N_3, v(a)+N_1+2N_3, v(a)+N_1+N_2+N_3, 
v(a)+v(\x)+N_3, \\
& v(a)+v(\ac\pnb-c)+N_1+N_3, v(\x)+N_1+N_3, 2N_1+N_2+N_3, 2v(\x), \\ &
2(N_1+N_2),  v(\x)+v(\ac\pnb-c)+N_1, v(\ac\pnb-c)+2N_1+N_2, \\ &
v(\x)+N_1+N_2, 2(N_1+N_3), v(b)+N_1+N_2+N_3,v(c)+2N_1+N_2, \\ &
v(c-\ac\pnb)+2N_1+N_3, v(a)+v(c)+N_1 +N_3,  v(c)+2N_1+N_3, \\ &
v(c)+v(\x)+N_1, v(c)+2N_1+v(\ac\pnb-c)
\bigr \}-w_1;
\endsplit
$$
and finally:
$$
\split
w_3 = \min \bigl\{ & 3v(a)+3N_3, v(a)+v(b)+N_1+N_2+N_3, 2v(a)+v(\x)+2N_3, \\ 
&
v(b)+2N_1+N_2+N_3, 3N_1+N_2+2N_3, 2v(a)+N_1+2N_3, \\ &
v(\x)+v(b)+N_1+N_2+N_3, v(\ac\pnb-c)+v(b)+2N_1+N_2+N_3, \\& 
v(\ac\pnb-c)+3N_1+N_2+N_3, v(a)+v(c)+2N_1+N_2+N_3, \\ &
v(a)+2v(\x) +N_3, v(c)+3N_1+N_2+N_3,
v(a)+v(\x) \\ &+N_1+N_3, v(a)+2N_1+N_2+N_3,
v(-\x^3+\aaa \x^2\pnab \\ & +\ab \x p^{2(N_1+N_2)}-c^2 p^{2N_1}\pnab +\ac c 
p^{2N_1}\pnab)
\bigr\} -w_1-w_2.
\endsplit
$$

We can similarly evaluate $\Omega_2, \dots, \Omega_6$ and get the following
values for these functions:
\proclaim {Proposition 2.3}
$$
\mu(\Omega_2) = \cases 0& \qquad \text {if }
M_2<V_1-(U_1+U_2)+N_1+N_2+N_3; \\
                            p^{V_1-(N_1+N_2+N_3)}& \quad \text{otherwise;}
                            \endcases
$$
where $V_1$ is defined to be
$$
V_1 :=\min \bigl \{N_1+N_3, v(a)+N_3, v(c)+N_1,v(\x),N_2+N_3 \bigr \}.
$$
\endproclaim
\proclaim{Proposition 2.4} $\Omega_3$ is $1$ if and only if
$$
\split & M_3 \geq N_1, M_3 \geq N_3,M_3+v(a)\geq N_1+N_2 ,  \\
& M_3 +v(\x) \geq N_1+N_2+N_3, M_3+v(c)\geq N_2+N_3
\endsplit
$$
and is $0$ otherwise.
\endproclaim

In an entirely similar fashion one can establish
\proclaim{Proposition 2.5 }
$\mu(\Omega_4)=  0 \text{ if } M_4<U_4+U_5-(W_1+W_2+W_3)+N_1+N_2+N_3$ and 
for all other values of $M_4$ we have
$$
\mu(\Omega_4)=p^{U_4+U_5-(N_1+N_2+N_3)};
$$
$$
\mu(\Omega_5) = \cases 0& \qquad \text{if }
M_5<V_4-(U_4+U_5)+N_1+N_2+N_3; \\
                            p^{V_4-(N_1+N_2+N_3)}& \quad \text{otherwise;}
\endcases
$$
and
$\Omega_6$ is $1$ if and only if
$$
\split & M_6 \geq N_1, M_6+v(a)\geq N_1+N_2 , M_6\geq N_2,  \\
& M_6 +v( - \y +\ac \pnab), M_6+v(\x) \geq N_1+N_2+N_3
\endsplit
$$
where the $U_i,V_i$ are defined as before.
\endproclaim

One can check that the only case that leads to non-monomial conditions on 
the entries of the matrix occurs when we evaluate $U_5:$
$$
\split
u_5=\min  \bigl \{ & 2v(a)+2N_3, v(a)+N_1+N_2+N_3, v(\x)+N_1+N_2,\\ &
v(a)+v(\x)+N_3, v(b)+N_1+N_2+N_3, 2N_1+N_2+N_3, \\&
v(\ac\pnb-c)+2N_1+N_2, 2v(\x)
\bigr \}-u_4.
\endsplit
$$

In order to show that evaluating the zeta function of the group depends on
counting points on the elliptic curve $\modp$ it will be
sufficient then to calculate the measure of the following set. For all
natural numbers $A,B,\X,C,F,G,H$ we need to find the
value of
$$
\split
\mu_{A,B,\X,C,F,G,H} :=\mu \bigl \{ &(a,b,c) \in \zp^3: 
v(a)=A,v(\x)=\X,v(c)=C, v(b)=B
\\
& v(\y+\ac \pnab)=G, v(\aaa\x+\ab\pnab)=H, \\&
v(-{\x}^3+ \aaa {\x}^2\pnab+ \ab \x p^{2(N_1+N_2)}- c^2p^{2N_1}\pnab \\&  
+\ac\y p^{2(N_1+N_2)})
=F \bigr \}.
\endsplit
$$
The first thing to notice is that as, writing $\Phi$ for the set $(b/c 
+p^{N_2-C}\zp)\cup p^A \zp^{\ast}$ we have
$$
\mu(\Phi) =\cases 0,&\text{if } B\not=A+C, N_2-C >\min\{A,B-C\}; \\
p^{-A}(1-p^{-1}), &\text{if } B\not=A+C, N_2-C \leq\min\{A,B-C\}; \\
p^{C-N_2}, &\text{if } B=A+C, A+C>N_2;\\
p^{-A}, &\text{if } B=A+C, A+C\leq N_2;
\endcases
$$
it is sufficient for us to evaluate
$$
\split
\mu_{\X,C,F,G,H} := \mu \bigl\{& (\x,c)\in \zp^2: v(\x)=\X, v(c)=C, 
v(\ac\pnab-c\pna)=G, \\ & v(\aaa\x+\ab\pnab)=H, v(-\x^3+\aaa\x^2\pnab +\ab\x 
p^{2(N_1+N_2)} \\& -c^2p^{2N_1} +\ac cp^{N_1}p^{2(N_1+N_2)})=F \bigr \}
\endsplit
$$
Thus by changing the value of $C,$ writing $(b,B)$ for $(\x,\X)$ and $N$ for 
$N_1+N_2$ we
need to calculate the value of
$$
\split
\mu_{B,C,F,G,H} := \mu \bigl \{ &(b,c) \in \zp^2: v(b)=B, v(c)=C,
v(-c +\ac p^N)=G, \\ & v(\aaa\x+\ab p^N)=H, v(-b^3+\aaa b^2p^N +\ab b 
p^{2N}\\& -c^2 p^N +\ac c p^{2N})=F\bigr \}
\endsplit
$$
We split the analysis into three sections:
\roster
\item $N\leq B,C;$
\item $B<N, B\leq C;$
\item $C<B,N.$
\endroster

\proclaim{Case 1} $N\leq B,C.$
\endproclaim

Setting $b^{\prime}=
b/{\pn}, c^{\prime}= c/{\pn}$ and replacing
$b^{\prime},c^{\prime}$ by $b,c$ respectively, we must evaluate in this 
instance
$$
\split
\mu \bigl \{& (b,c) \in \zp^2: v(b)=B, v(c)=C, v(\ac -c)
= G,
\\& v(\aaa b+\ab)=H, v(b^3-\aaa b^2 -\ab b
+c^2  -\ac c)= F  \bigr \}
\endsplit
$$
Notice that the calculations that follow assume that $\aaa\ab\ac \not= 0.$ 
The special cases that follow when this is not the case are all handled in 
the same way and so we suppress the details.

Notice also that the above set can be expressed as a Boolean combination of 
sets of the form
$$
\split
\bigl \{& (b,c) \in \zp^2: v(b)=B, v(c)=C, v(\ac -c)
\geq G,
\\& v(\aaa b+\ab)\geq H, v(b^3-\aaa b^2 -\ab b
+c^2  -\ac c)\geq F  \bigr \}.
\endsplit
$$
We write $d(B,C,F,G,H)$ for the measure of this set and evaluate this. The 
first thing to notice is:

\proclaim{Lemma 2.6}
Suppose $B,C>0.$ Then $G=0=H$ otherwise $d(B,C,F,G,H)$ is $0.$ When $G=0=H$ 
then
the measure depends upon the value of $F$ in
relation to that of $B$ and $C.$ Namely:
\roster
\item  if $F \leq \min \{B,C\}$ then
$d(B,C,F,0)=p^{-C}p^{-B}(1-p^{-1})^2;$
\item if $ F<\min
\{ B,C\} $ then $d(B,C,F,0)=0$ unless $B=C$ when
$d(B,C,F,0)=p^{-F}p^{-C}(1-p^{-1}).$
\endroster
\endproclaim

Next consider what happens for $C>0,B=0.$ As in Lemma 2.6 we require $G=0$ 
in order to get a non-zero value for the measure. We encapsulate what 
happens in this instance in the next

\proclaim{Lemma 2.7}
Suppose that $C>0=B=G.$ then
\roster
\item $H>0 \implies F=0$ or the measure is zero. In the case $H>0,F=0$ we 
must evaluate
$$
\mu\bigl \{(b,c)\in \zp^2:v(b)=0, v(c)=C, v(\aaa b+\ab)\geq H \bigr\}
$$
which is dependent upon the number of points on the line $\{\aaa b+\ab 
=0\}(\Bbb F_p).$ This is uniform in $p$ and so can be neglected$;$
\item If $G=0=H$ then $1 \leq F
\leq C \Rightarrow d(0,C,F,0,0)=2p^{-F}p^{-C}(1-p^{-1}).$
\endroster
\endproclaim

Next we consider the case $B>0, C=0.$ As previously it is immediate that to 
get a non-zero value for the measure we require that $H=0.$ When $H=0$ we 
evaluate
$$
\split
d^{\prime}(B,0,F,G,0)= \mu \bigl\{&v(b)=B,v(c)=0, v(\ac-c)= G, \\& 
v(b^3-\aaa b^2-\ab b +c^2-\ac c)\geq F \bigr\}.
\endsplit
$$

\proclaim{Lemma 2.8}
\roster
\item If $F>\min \{B,G\}$ then $d^{\prime}(B,0,F,G,0)= 0;$
\item If $F\leq \min \{B,G\}$ then $d^{\prime}(B,0,F,G,0)$ depends upon the 
number of points on the line $\{\ac -c =0\}(\Bbb F_p)$ and so is uniform in 
$p.$
\endroster
\endproclaim

Finally we must consider what happens when $B=C=0.$ In this case we want to 
calculate a value for
$$
\split
d_{0,0,F,G,H}&=\mu\bigl \{(b,c) \in \zp^2: v(b)=0, v(c)=0, v(\ac -c) \geq G,
\\
&v(\aaa b+\ab)\geq H, v(b^3+\aaa b^2 +\ab b
+c^2  +\ac c)\geq F  \bigr \}
\endsplit
$$
The dependence on the varieties described in Theorem 1 will be born out of 
the following

\proclaim{Lemma 2.9}
Let $K \geq 1,$ let $(b,c) \in (\frac {\Bbb Z}{p^K\Bbb Z})$ with
$$
A_1b^3+A_2b^2+A_3b+A_5c^2+A_6c \equiv   0 \mod  p^K
$$
where here $p$ is a prime dividing neither  the discriminant nor the 
coefficients of the curve. Then
there exist $p$ pairs $(b_1,c_1) \in
(\frac {\Bbb Z}{p^{K+1}\Bbb Z})$ so that $b \equiv   b_1,\, c\equiv   c_1
\mod p^{K+1}$ and
$$
A_1b_1^3+A_2b_1^2+A_3b_1+A_5c_1^2+A_6c_1 \equiv   0 \mod  p^{K+1}.
\tag 2.2
$$

\endproclaim

\demo{Proof}
Setting $b_1=b+\beta p^K$ and $c_1=c+\gamma p^K$ we want to count pairs
$(\beta, \gamma ) \in \{0,\dots ,p-1\}^2$ so that (2.2) is satisfied.
Expand
this equation and notice that we may write
$$
A_1b^3 +\dots +A_6c = tp^K
$$
for some $t \in \Bbb N,$ and then it follows that we are looking for
solutions of the linear congruence (in terms of $\beta$ and $\gamma$)
$$
t+\beta (3b^2A_1+2bA_2+A_3) +\gamma (2A_5+A_6) \equiv   0  \modp
$$
The only way this congruence cannot have $p$ solutions is when both the
coefficients of $\beta $ and $\gamma$ are zero $\modp.$ But this
happens only when $p$ divides the discriminant of the curve contradicting
the hypothesis we made. Thus the result is proved.

\enddemo
We will split the calculation of the measure into a case analysis dependent 
upon the values of $F,G,H:$
\roster
\item $F,G,H =0;$
\item $G,H=0;$
\item $F,H=0;$
\item $F,G=0;$
\item $F=0;$
\item $G=0;$
\item $H=0;$
\item $F,G,H \not=0.$
\endroster
\
1) We have $d(0,0,0,0,0) = \mu (\zp^\ast \times \zp^\ast)$ which is uniform 
in $p$ and so can be neglected.
\newline
2) It follows simply from Lemma 2.9 that
$$
d(0,0,F,0,0) =p^{-F+1}d(0,0,1,0,0)= p^{-F+1}(|E(\Bbb F_p|-1)
$$
\newline
3) We have $d(0,0,0,G,0) =(1-p^{-1})p^{-G}.$ Notice that this expression 
actually involves counting points on the line $\{\ac-c=0\}(\Bbb F_p)$ but 
this is suppressed due to the uniformity of this variety.
\newline
4) In an identical fashion, we have $d(0,0,0,0,H)=(1-p^{-1})p^{-H}.$
\newline
5) $d(0,0,0,G,H) = p^{-G-H};$ here we are counting points on the 
intersection $\{\ac-c=0\}\cap \{\aaa b+\ab=0\}(\Bbb F_p).$
\newline
6) $d(0,0,F,0,1) = p^{-F+1}\bigl (E\cap \{\aaa b+\ab =0\}(\Bbb F_p)\bigr )$ 
and the case for a general $H$ follows as a simple recurrence relation.
\newline
7) $d(0,0,F,1,0) = p^{-F+1}\bigl (E\cap \{\ac-c =0\}(\Bbb F_p)\bigr )$ and 
the case for a general $G$ follows as a simple recurrence relation.
\newline
8) $d(0,0,F,1,1) = p^{-F+1}\bigl (E\cap \{\aaa b+\ab =0\}\cap \{\ac-c 
=0\}(\Bbb F_p)\bigr )$ and the case for general $G,H$ follows as a simple 
recurrence relation.
\newline

This completes the case $N\leq B,C.$ Notice that this is as stipulated
by the work of du Sautoy and Grunewald in \cite{duSG1} in that
finitely many varieties, and their intersections, arise in the evaluation of
the local zeta function.

\proclaim{Case 2} $B<N, B\leq C.$
\endproclaim

Setting $b^{\prime}=
{\pn}/b, \, c^{\prime} = c/b$ and replacing as before we must
evaluate, for $B\geq1, C\geq 0$
$$
\split
\mu  \bigl \{ (b,c) & \in  \zp^2 :v(b)=B,v(c)=C, v(b\ac -c) =
G, \\& v(\aaa+\ab b)=H,
v(1 +\aaa b +\ab  b^2 +c^2b+\ac cb^2)= F \bigr \}.
\endsplit
$$
Using the same notation as before we can see that $F,H>0 \implies 
d(B,C,F,G,H)=0.$ The only non-trivial measures arising here are contained in
\proclaim{Lemma 2.10} Suppose that $F=0=H.$ Then
\roster
\item $G>\min \{B,C\}\implies  d(B,C,0,G,0) =0$
\item $\min \{B,C\} \geq G \implies p^{-B}p^{-C}(1-p^{-1})^{-2}$
\endroster
\endproclaim

\proclaim{Case 3} $C<B,N.$
\endproclaim

As is now becoming familiar, we set $b^{\prime}=
b/c , \, c^{\prime}=  {\pn}/c$ and relabelling as before we
must evaluate for $B,C>0$
$$
\split
\mu \bigl\{  (b,c)& \in \zp^2:v(b)=B,v(c)=C, v(\ac c-1)=G,\\& v(\aaa b+\ab 
c)=H,  v(b^3-\aaa b^2c-\ab bc^2 +c-\ac c^2)=F \bigr\}.
\endsplit
$$
We again immediately notice that $G>0$ will give us a set of measure $0.$
The remaining cases, with the usual notation, are encapsulated in the 
following

\proclaim {Lemma 2.11}
\roster
\item If $F \leq \min \{3B,C\}$ and $H\leq \min\{B,C\}$ then
$$
d(B,C,F,0,H) = p^{-C}p^{-B}(1-p^{-1})^2;
$$
\item if either $F>\min\{3B,C\}$ and $3B\not=C,$ or $H>\min\{B,C\}$ and 
$B\not=C$ then $d(B,C,F,0,H)=0;$
\item if $F>\min \{3B,C\}, 3B=C$ and $H\leq \min\{B,C\}$ then
$$
d(B,3B,F,0,H) =p^{-F}p^{-B}(1-p^{-1});
$$
\item if $H>\min\{B,C\}, B=C$ and $F\leq \min\{3B,C\}$ then
$$
d(B,C,F,0,H)=p^{-H}p^{-B}(1-p^{-1}).
$$
\endroster
\endproclaim

So to finish I will briefly outline why the calculations I have performed
lead to the theorem stated in the Introduction. This is exactly
as contained in \cite{duS3} and so I only include it for completeness. The
calculations carried out here are sufficient to prove

\proclaim{Proposition 2.12}
There exists a finite partition $\cup_{i \in S} \Delta_i$ of $\Bbb R^9$
defined by linear inequalities with coefficients in $\Bbb Q$ and
for all $i \in S$ polynomials $P_i,Q_i,R_i,S_i,T_i $ in $\Bbb Q(X)$ and 
linear
functions $\alpha_i,\beta_i,\gamma_i,\delta_i,\epsilon_i$ so that if
$$
\Delta = (A,B,\X,C,F,G,H,N_1,N_2) \in \Bbb N^9 \cap \Delta_i
$$
then
$$
\split
\mu \{ \dots \} = P_i(p)p^{\alpha_i(\Delta)} &+ Q_i(p)|E(\Bbb
F_p)|p^{\beta_i(\Delta)} + R_i(p)|E \cap M_1 (\Bbb F_p)|p^{\gamma_i(\Delta)}
\\
&+
S_i(p)|E\cap M_2 (\Bbb F_p)|p^{\delta_i(\Delta)}+T_i(p)|E\cap M_1\cap 
M_2(\Bbb F_p)|p^{\epsilon_i(\Delta)}.
\endsplit
$$
\endproclaim

{}From this result, together with the values we worked out for the functions
$\Omega _1, \dots , \Omega_6$ one can deduce that there exists
a finite partition $\cup_{i \in S} \Delta_i$ of $\Bbb R^{18}$ defined by
linear inequalities with coefficients in $\Bbb Q$ and
for all $i \in S$ polynomials $P_i,Q_i,R_i,S_i $ in $ \Bbb Q(X)$ and linear
functions 
$\alpha_i,\beta_i,\gamma_i,\delta_i,\epsilon_i,a_i,b_i,c_i,d_i,e_i$ so that 
if
$\Lambda
=(A,B,\X,C,F,G,H,M_1,\dots ,M_6,N_1,N_2,N_3)$ then
$$
\split
\zlp =&\sum_{i \in S} \sum_{\Lambda \in \Bbb N^{18}\cap \Delta_i}
P_i(p)p^{\alpha_i(\Lambda)+a_i(\Lambda)s} + Q_i(p)|E(\Bbb
F_p)|p^{\beta_i(\Lambda)+b_i(\Lambda)s} +  \\ & R_i(p)
|E \cap M_1 (\Bbb F_p)|p^{\gamma_i(\Lambda)+c_i(\Lambda)s} +
S_i(p)|E\cap M_2 (\Bbb F_p)|p^{\delta_i(\Lambda)+d_i(\Lambda)s} \\&
T_i(p)|E\cap M_1\cap M_2 (\Bbb F_p)|p^{\epsilon_i(\Lambda)+e_i(\Lambda)s}.
\endsplit
$$

Adding all this together is almost sufficient to prove Theorem 1.
It is merely necessary to check that the rational function
$P_2$ is not identically equal to $0.$

\vskip 0.5in
\heading {\bf {Section 3: Completion of the proof of Theorem 1}}
\endheading
\vskip 0.25in

Although the Theorem is now complete, in order to show that this collection
of nilpotent groups really does encode the arithmetic of the
elliptic curves it is necessary to show that the rational function we have
called $P_2$ is non-zero. As things stand we have merely shown
the existence of such a function without saying anything about what it looks
like. To show the function is non-zero, it is sufficient to
show that counting subalgebras of some small $p$-power index in $L$ is
dependent on counting points on the reduction of the elliptic
curves in question. This is a simple exercise in solving some congruences
$\modp$.

So again let $L=L(E)$ be the Lie algebra with presentation as described
earlier. Throughout this calculation we will assume we are dealing
with a prime $p$ not dividing $\aaa\ab\ac$ as this will simplify
greatly the work involved. Then to count ideals of
index $p^5$ say, it will be sufficient to count the number of pairs of
matrices $((m_{ij}), \pmatrix n_1 & a & b \\ 0 &
n_2 & c \\ 0 & 0 & n_3 \endpmatrix )$ so that the following four conditions
are satisfied:

\roster
\item $ m_{ij},n_k \in \Bbb Z;$
\item $ 0\leq m_{ij}<m_{jj}, \quad 0\leq a< n_2,\quad 0\leq b,c<n_3;$
\item $ m_{jj}=p^{M_j}, n_j=p^{N_j}, \quad M_1+\dotsb N_3 =5 ;$
\item for $ i=1, \dots , 6;$
$$
j=1,2,3 \Rightarrow (m_{i4},m_{i5},m_{i6})C(j)N^+ =
(\alpha_0,\beta_0,\gamma_0)n_1n_2n_3;
$$
$$
j=4,5,6 \Rightarrow (m_{i1},m_{i2},m_{i3})C(j)N^+ =
(\alpha_0,\beta_0,\gamma_0)n_1n_2n_3.
$$
\endroster

Now the first thing to notice is that if we work out the left hand side of
condition (4) for all the relevant values of $i,j$ and $\alpha_i$
then we can immediately deduce the following

\proclaim{Lemma 3.1}
$M_1,\dots, M_6 \geq N_1; M_3,M_6 \geq N_2; M_2,M_5 \geq N_3.$ Hence $N_1=0$
and $0\leq N_2,N_3 \leq 1$ and furthermore $N_2+N_3=0,1.$
\endproclaim

To make the analysis that follows more tractable, we again split the working 
into several separate cases.
\proclaim {Case 1}$N_2=N_3=0$
\endproclaim
This is easily dealt with. In this case $N=\text {
Id}_3$, no conditions arise from (4) and we merely have to count all
matrices $(m_{ij})$ that can occur. This is uniform and polynomial in $p$
and so can be encompassed under the umbrella of a rational
function of $p,p^{-s}$ in the evaluation of the zeta function. As such it
doesn\rq t concern us here.

\proclaim{Case 2}$N_2=1, N_3=0.$
\endproclaim

It follows from the conditions stipulated above that $N_3=b=c=0$ and $0\leq
a \leq p-1.$ Then condition (4) becomes
$$
(m_{i4},m_{i5},m_{i6}) \pmatrix \aaa p & -a\aaa & \ab p \\ p & -a & 0 \\ 0 &
1 & \ac p \endpmatrix \equiv   0   \modp ;
$$

$$
(m_{i1},m_{i2},m_{i3}) \pmatrix \aaa p & -a\aaa & ap \\ \ac p & -a & 0 \\
0 & 1 & 0 \endpmatrix  \equiv   0  \modp ;
$$

$$
(m_{i4},m_{i5},m_{i6}) \pmatrix p & -a & 0 \\ 0 & 0 & p \\ 0 & 0 & 0
\endpmatrix \equiv   0  \modp ;
$$

$$
(m_{i4},m_{i5},m_{i6}) \pmatrix 0 & 1 & 0  \\ 0 & 0 & 0 \\ p & -a & 0
\endpmatrix  \equiv   0  \modp ;
$$

$$
(m_{i1},m_{i2},m_{i3}) \pmatrix p & -a & 0 \\ 0 & 0 & p \\ 0 & 0 & 0
\endpmatrix \equiv   0  \modp ;
$$

$$
(m_{i1},m_{i2},m_{i3}) \pmatrix 0 & 1 & \ac p \\ 0 & 0 & 0 \\ p & -a & 0
\endpmatrix  \equiv   0  \modp .
$$

So if $a=0$ then the only thing these matrices tell us is that $m_{i6}
\equiv   0 \modp$ from whence it follows that
$m_{i4}=m_{i3}=m_{i1} \equiv   0 \modp$ and thus the number of matrices in
this case is
again uniform and polynomial in $p.$ If on the other hand $1\leq a \leq p-1$
then we see that $m_{i1} \equiv   0   \equiv   m_{i4}
\modp$ and so $\underline a = (1,0,1,1,0,1)$ and again we will get a uniform
expression.

\proclaim{Case 3} Suppose finally that $N_3=1, \, N_2=0=a.$
\endproclaim

In this case the
conditions the matrices must satisfy become
$$
\align
(\ab - b\aaa )m_{i4}+bm_{i5}+(\ac  - c )m_{i6}&  \equiv   0  \modp ; \tag
3.1
\\
bm_{i4}+m_{i5}& \equiv    0  \modp ;   \tag 3.2\\
-cm_{i4}+ bm_{i6}& \equiv   0  \modp ;  \tag 3.3\\
(\ab+b\aaa)m_{i1}+b m_{i2}-cm_{i3}& \equiv   0  \modp ;  \tag 3.4 \\
bm_{i1}+m_{i2}& \equiv   0  \modp ;   \tag 3.5\\
(- c+\ac b)m_{i1}+bm_{i3}& \equiv   0  \modp . \tag 3.6
\endalign
$$

Recall that we already know that $p$ divides $m_{22},\, m_{55}$ by Lemma
3.1 and this is confirmed by equations (3.2), (3.5).
Setting $i=3,6$ allows us to deduce the following congruences involving
entries on the diagonal of $(m_{ij});$
$$
\align
-bm_{66}  & \equiv   0  \modp ; \tag 3.7\\
(\ac  -c )m_{66} & \equiv   0  \modp ; \tag 3.8\\
-bm_{33} & \equiv   0  \modp ; \tag 3.9\\
cm_{33} & \equiv   0  \modp . \tag 3.10
\endalign
$$
So we do a case-by-case analysis, dependent on the values of $ b,
\ac-c
,c  \modp.$

\proclaim {Subcase 3.2} $ b \not \equiv   0  \modp$ so that $p$ does
not divide $b.$
\endproclaim
We can immediately see that $p$ divides
$m_{33},m_{66}$ and so $\underline a =(0,1,1,0,1,1).$ Setting $i=1$ tells us
that
$$
\aligned
m_{12} \equiv  -b  \modp ; \qquad \qquad \\
m_{13} \equiv  \frac {-( c-\ac)}{ b}  \modp ; \qquad \qquad \\
(\ab +\aaa b )+ b m_{12}-cm_{13} \equiv  0  \modp  \qquad \qquad
\endaligned
$$
and so we have uniquely determined $m_{12},m_{13}$ in terms of $b,c$ from
which it follows by
elementary analysis that $(b,c) \in E(\Bbb F_p).$

An entirely similar process gives us the same information about
$m_{45},m_{46}.$ So now to deduce that in this case we have no alternative
but to count points on the reduction of the elliptic curve it will be
sufficient to demonstrate what values the other entries in the
matrix $M$ can take.

Recall we know from the four conditions that $m_{i4}=0\, \forall i<4.$ So we
need to investigate how to determine the remaining values of
$m_{i2},m_{i3},m_{i5},m_{i6}.$ This is easily done simply by examining the
equations for values of $i$ running from $1$ through $5$ and
one can see that all the outstanding values are uniquely determined by the
choice of point $(b,c)$ on the reduction of the curve.

\proclaim {Subcase 3.3} $b \equiv  0  \modp.$
\endproclaim
The $6$ conditions then become
$$
\aligned
\ab m_{i4}+(\ac - c)m_{i6}& \equiv  0  \modp; \\
m_{i5}& \equiv  0  \modp; \\
-cm_{i4}& \equiv  0  \modp; \\
\ab m_{i1}-cm_{i3}& \equiv  0  \modp; \\
m_{i2}& \equiv  0  \modp; \\
(\ac -  c)m_{i1}& \equiv  0  \modp.
\endaligned
$$

{}From these equations it is possible to show that in  two of the cases that
can occur, namely 1) $c,\ac - c \not \equiv  0  \modp$ and 2) $c
\equiv  0  \modp , \ac - c  \not \equiv  0  \modp$ that the fact that
we are working with ideals of index $p^5$ means that no such
matrices can occur. It seems reasonable however that if we increase the
exponent of $p$ then matrices will occur that bear witness to
these congruences. However suppose we look at the final possible case $c
\not \equiv  0  \modp, \ac - c \equiv  0  \modp.$ In this
case the conditions become
$$
\aligned
\ab m_{i4}& \equiv  0  \modp; \\
m_{i5}& \equiv  0  \modp; \\
cm_{i4}& \equiv  0  \modp; \\
\ab m_{i1}-cm_{i3}& \equiv  0  \modp; \\
m_{i2}& \equiv  0  \modp
\endaligned
$$
and so the only constraints the coefficients are bound by are
$$
m_{22},m_{33},m_{44},m_{55} \equiv  0  \modp.
$$
It can again be checked
that this leads to a polynomial uniform expression in $p.$ This finishes the
calculation and demonstrates that it really is necessary to count points on
the elliptic curve, there is no quirk which ensures a simple expression
after all. Thus Theorem 1 is proved.

\vskip 0.5in
\heading {\bf {Section 4: Proof of Theorem 2}}
\endheading
\vskip 0.25in
First let us recall what Theorem 2 stated:
\proclaim{Theorem 2}
Let $C$ be a curve of genus 2 over $\Bbb Q$ with a rational point in $\Bbb
Q.$ Then there exists a 15 generated, class 2 torsion-free nilpotent group G
and associated Lie algebra L, $l \in \Bbb N,$ rational functions $P,
Q_1,\dots, Q_l \in \Bbb Q(X,Y)$ and varieties $V_1, \dots, V_l$ so that for
almost all primes $p$
$$
\zeta_{G,p}^{\triangleleft}(s) =\zeta_{L,p}^{\triangleleft} (s) = |C(\Bbb
F_p)|P(p,p^{-s}) + \sum_{i=1}^l |V_i(\Bbb F_p)| Q_i(p,p^{-s}).
$$
Furthermore it is strictly necessary to count points on the reduction of the 
curve $C$ in the evaluation of the zeta function.  In particular, the $V_i$ 
are varieties of genus smaller than 2 and $C$ occurs in the subring of the 
motivic zeta function one can associate to the group $G.$
\endproclaim
As this result follows along very similar lines to Theorem 1, a lot of
details will be swept under the carpet. It is known (see \cite {M2, Ch1} for
example) that every curve of genus 2 over $\Bbb Q$ is of the form
$$
Y^2=a_0X^6+ \dots +a_6
$$
and so every curve of genus 2 with a rational point is of the form
$$
Y^2 + bY=a_0X^6+\dots + a_5X
$$
or projectively of the form
$$
Y^2Z^4+bYZ^5=a_0X^6+\dots + a_5XZ^5
$$
and so is expressible as the determinant of
$$
G:= \pmatrix Y & X & \beta_1X & 0 & \beta_2X+\beta_3Z & \beta_4Z \\ 0 & Z &X
& \beta_5Z & 0 & 0 \\ 0 & 0 & Z & X & 0 & 0 \\ 0 & 0 & 0 & Z & X & 0 \\ 0& 0
& 0 & 0 & Z & X \\ \beta_6X & 0 & 0 & 0 & 0 & Y+\beta_7Z \endpmatrix .
$$
Thus we will count ideals of $p$-power index in
$$
L:= \langle A_1, \dots, A_6,B_1, \dots,
B_6,X,Y,Z:(A_i,B_j)=g_{ij}(X,Y,Z)\rangle .
$$
For an algebra of this size, it becomes very difficult to evaluate the
integral that would give us a full description of the zeta function encoding
the ideal structure of $L.$ Thus for this algebra, we will content ourselves
with evaluating the coefficients of the zeta function for small powers of
$p$ in order to demonstrate that evaluating the zeta function does depend
upon the number of points on the genus 2 curve as claimed. Given the work of
du Sautoy and Grunewald in \cite{duSG1} this will be sufficient to prove
Theorem 2. This will follow exactly as for the elliptic curve example
already considered, and in fact the details are very similar also. Recall
that in order to count ideals of index $p^n$ in $L$ it is necessary and
sufficient to count all pairs of matrices $((m_{ij}), \pmatrix n_1 & a & b\\
0 & n_2 & c \\ 0 & 0 & n_3 \endpmatrix )$ so that the following four
conditions are satisfied:
\roster
\item $ m_{ij},n_k \in \Bbb Z;$
\item $ 0\leq m_{ij}<m_{jj}, \quad 0\leq a< n_2,\quad 0\leq b,c<n_3;$
\item $ m_{jj}=p^{M_j}, n_j=p^{N_j}, \quad M_1+\dotsb +N_3 =n ;$
\item for $ i=1, \dots , 6;$
$$
j=1,\dots ,6 \Rightarrow (m_{i7}, \dots ,m_{i12})C(j)N^+ =
(\alpha_0,\beta_0,\gamma_0)n_1n_2n_3;
$$
$$
j=7, \dots ,12 \Rightarrow (m_{i1}, \dots ,m_{i6},)C(j)N^+
=(\alpha_0,\beta_0,\gamma_0)n_1n_2n_3.
$$
\endroster
where the matrices $C(1), \dots ,C(12)$ are defined as in the elliptic curve
example. Then if $a_n$ denotes the number of ideals of index $p^n$ in the
algebra $L,$ then
$$
a_n=\sum_{<(\underline a, \underline b)>=n}c_{(\underline a, \underline
b)}p^{12(b_1+b_2+b_3)}
$$
where $c_{(\underline a, \underline b)}$ denotes the number of pairs of
matrices satisfying the above conditions with diagonal entries
$p^{a_i},p^{b_i}$ respectively. Also note that by $<(\underline a,\underline
b)>$ we mean the sum of the entries in the vectors. Again we immediately get
some restrictions on the values the diagonal entries of the matrices can
take which are included in the next
\proclaim {Lemma 4.1}
\roster
\item $a_1, \dots ,a_7, a_9, \dots ,a_{12} \geq b_1;$
\item $a_2,a_4,a_5, a_{12} \geq b_3;$
\item $a_6,a_{12} \geq b_2.$
\endroster
\endproclaim
We consider the case $n=11$ as in this instance it is a simple task to
demonstrate the dependence on the curve of the number of ideals of given
index. The first thing to notice is that the above Lemma forces $b_1=0.$ We
consider the case
$$
(\underline a, \underline b)>= (0,1,1,1,1,1,0,1,1,1,1,1,0,0,1).
$$
In other words, we are counting pairs of matrices $(M,N)$ in which
$$
N=\pmatrix 1 & 0 & b \\ 0 & 1 & c \\ 0 & 0 & p \endpmatrix ,
$$
the diagonal entries of $M$ are $(1,p,p,p,p,p,1,p,p,p,p,p)$ and condition
$(4)$ above has now reduced to (replacing $b$ by $-b$)
$$
\align–c m_{i7} +b m_{i8} +\beta_1b m_{i9} + (\beta_2b +\beta_3) m_{i11}
+\beta_4 m_{i12}& \equiv 0  \modp ;\\
m_{i8} +bm_{i9} +\beta_5m_{i10} &\equiv 0  \modp ; \\
m_{i9} +bm_{i10} &\equiv 0  \modp ; \\
m_{i10} +bm_{i11} &\equiv 0  \modp ;\\
m_{i11} +bm_{i12} &\equiv 0  \modp ;\\
\beta_6bm_{i7}+ (\beta_7-c)m_{i12} &\equiv 0  \modp ;\\
-cm_{i1} +bm_{i6} &\equiv 0  \modp ;\\
bm_{i1} +m_{i2}&\equiv 0  \modp ;\\
\beta_1bm_{i1}+b m_{i2}+ m_{i3}&\equiv 0  \modp ; \\
\beta_5m_{i2}+ bm_{i3} +m_{i4}&\equiv 0  \modp ; \\
(\beta_2b+\beta_3)m_{i1}
+bm_{i4} +m_{i5}&\equiv 0  \modp ;\\
\beta_4m_{i1} +bm_{i5} +(\beta_7-c)m_{i6}&\equiv 0  \modp .\endalign
$$
Notice that if $b,c=0$ then these conditions reduce even further and we get
a uniform number of pairs of matrices, regardless of the prime $p.$ The
interesting case is when $b$ and $c$ are non-zero. In this case, as for the
case of the elliptic curve example, it is simple to see that the number of
pairs of matrices that occur is dependent on counting points on the
reduction of the curve. For $M$ is uniquely determined and the above
conditions imply that the number of various $N$ that can occur is $|C(\Bbb
F_p)|-1.$ Theorem 2 now follows as in the work of du Sautoy
in \cite{duS4} which we mirrored in Section 3 of this paper.

{\bf Remark.} It is possible to complete this working and get a full 
description of $a_{11}.$ However for the purposes of the Theorem we have 
sufficient detail.

\vskip 0.5in
\heading{\bf {Section 5: Counting All Subalgebras}}
\endheading
\vskip 0.25in

We have given a fairly complete description of the zeta
function counting ideals of the Lie algebras considered in the proofs of
Theorems 1 and 2. We now digress slightly and consider the problem of
counting all subalgebras in the Lie algebra $L$ associated to the elliptic
curve of Theorem 1. It is natural, given the presentation of
$L,$ to suspect that evaluating the zeta function counting all subalgebras
of $L$ will depend on counting points on the same varieties. Again it is
known, from \cite{dSG1}, that

$$
\zeta_{L,p}^{\leq}(s) = (1-p^{-1})^{-9}\int_{W_p} |m_{11}|^{s-1} \dotsb
|m_{99}|^{s-9} |dx|
$$
and we again simplify to write as

$$
\zeta_{L,p}^{\leq}(s) = (1-p^{-1})^{-9}\int_{W_p} |m_{11}|^{s-1} \dotsb
|n_{3}|^{s-9} |dm|\cdot |dn|
$$
where now $W_p$ consists of upper triangular matrices $M \in Tr_6(\zp)$ so
that for all $1 \leq i,j \leq 6$,
$$
\underline m_i (\sum_{l=j}^6 m_{jl}D(l))N^+ \in n_1n_2n_3 \zp^3
$$
where
$$
D(1) = \pmatrix 0 & 0 & 0 \\ 0 & 0 & 0 \\ 0 & 0 & 0 \\ -\aaa & 0 & -\ab \\
-1 & 0 & 0 \\ 0 & -1 & -\ac \endpmatrix
$$
and $D(2), \dots ,D(6)$ are defined similarly. As we have done above, it is
possible to explicitly count all subalgebras of $L$ of small $p$-power
index, and these calculations, which are too lengthy to be included here,
head me to pose a

\proclaim {Problem 5.1} Given the algebra $L(E)$ associated to the
elliptic
curve $E$ and nilpotent group $G(E)$ do there exist rational functions $P_1, 
P_2,P_3,P_4,P_5$ so that for almost
all primes
$$
\split
\zeta_{G,p}^{\leq} (s) = \zeta_{L,p}^{\leq} (s) & = P_1(p,p^{-s}) +|E(\Bbb 
F_p)| P_2(p,p^{-s})+|M_1\cap
E(\Bbb F_p)|P_3(p,p^{-s}) \\ & +P_4(p,p^{-s})|M_2 \cap E(\Bbb F_p)| 
+P_5(p,p^{-s})|M_1\cap M_2 \cap E(\Bbb F_p)|
\endsplit
$$
for the same lines $M_1,M_2?$
\endproclaim

The stumbling block to proving this isn\rq t one of conception: all the
machinery would appear to be in place. However the conditions lead to a very
complicated case analysis which I have not yet carried out.

This does lead one to another:

\proclaim{Question 5.2}
Given a finite dimensional Lie ring $L$ over $\bbbz$ do the same varieties 
always
arise when one evaluates either the local zeta function counting all prime 
power index subalgebras or the local zeta function only counting all prime 
power index ideals?
\endproclaim

A proof would probably come from
understanding the associated motivic zeta function better.

\vskip 0.5in
\heading {\bf {Section 6: Questions arising from this work}}
\endheading
\vskip 0.25in

We have extended du Sautoy\rq s work to produce a larger class of curves
whose arithmetic is encoded in the subgroup structure of some nilpotent
groups. Although in the genus 2 example we already see that things rapidly
become more complicated as the degree of the curve increases, we may still
in theory ask to what extent this method of producing curves as determinants
holds good.

In 1921, Dickson \cite{Di} considered the problem over $\Bbb C$ and gave a
description of all homogeneous polynomials arising as the determinant of a
matrix with linear entries; his methods were somewhat ad hoc. More recently,
Beauville \cite{Be} has used the theory of Cohen-Macauley sheaves to show in
fact that any curve over $\Bbb Q$ can be written as the determinant of a
matrix of linear forms. Thus it is in theory possible to extend further my
examples and produce any curve as a determinant.

Here we have considered curves of small genus as they are classes of curves
with a nice general description. It may be that one can define other classes
of curves with general equations of this type, but it known for instance
\cite {M2} that when one considers curves of genus 3, such curves do not
have a description of a similar kind. There is no general formula giving
every such curve.

Also notice we have stipulated that the curves we consider have a rational
point. This is to keep notation as simple as possible. For example, it does
not appear possible to write down a determinant giving an arbitrary elliptic
curve; it seems that as one considers curves, one must consider alternative
styles of presentation. The class of curves with a rational point contains a
large proportion of all elliptic curves, conjecturally 70\% of them \cite 
{BM,W}, and
this class does have a nice expression as a determinant as we have seen.

Similar work to that contained here has been carried out by Christopher Voll
\cite {V} who has also considered the problem of constructing groups whose
subgroup structure encodes information about the reduction of some plane
curves. He considers more generally curves over an algebraic number field,
whose representation as a determinant is well known \cite {Di}. He is able
to give expressions for the zeta function of a Lie algebra whose Lie
structure is defined similarly to that contained here. As such he
demonstrates a relationship between plane curves over a number field and
zeta functions counting certain restricted types of subalgebras of a Lie
algebra defined over $\Bbb Q.$ For more details consult his thesis.

Voll also completes the calculation for the zeta function counting points on 
the elliptic curve $E:Y^2=X^3-X.$ By explicitly evaluating the rational 
functions, and applying a functional equation for the number of points on 
$E(\Bbb F_p)$ he is able to demonstrate the existence of a functional 
equation for this zeta function. One can ask whether this phenomenon will 
hold in full generality, for instance for the zeta functions considered 
here. Denef and Meuser \cite {DM} have shown that the Igusa zeta function 
has a functional equation; this zeta function is a special case of du Sautoy 
and Grunewald\rq s cone integral with an empty cone condition. It is 
possible to construct a cone condition so that the associated cone integral 
does not satisfy a  functional equation, but can such a cone condition come 
from a presentation for a nilpotent group? Or do the cone conditions arising 
from group presentations all have the necessary symmetry to ensure the 
existence of a functional equation for all group zeta functions? I thank 
Marcus du Sautoy for suggesting this reasoning to me.

To end this paper, I will now note the determinants arising from the 
$3\times 3$ minors of the matrix $(S_1,S_2)$ in the calculation of the zeta 
functions in Section 1 of this paper.

\vskip 0.5in
\heading {\bf {Appendix}}
\endheading
\vskip 0.25in

In the interests of completeness, we include here the determinants arising 
from the
$3\times 3$ minors of the matrix $(S_1,S_2).$ By repeatedly applying the 
condition $\min \{v(X+Y),v(X)\} =\min \{v(X),v(Y)\}$ and noticing that we 
have stipulated that $p$ does not divide $\aaa\ab\ac$ it is relatively 
straightforward to see that the value
for $W_3$ is as contained in the main body of the text. The same process 
enables one to evaluate $W_2,U_2,U_5$ but in the interests of brevity we 
suppress the details. Throughout $(a_1,a_2,a_3)$ will
denote the determinant arising from the matrix formed from the $a_1,a_2,a_3$ 
columns of
$(S_1,S_2)$ for $a_i \in \{1,\dots, 6\}:$

\roster
\item"{$(1,2,3)$}" $a^3p^{3N_3}$
\item"{$(1,2,4)$}" $abp^{N_1+N_2+2N_3}$
\item"{$(1,2,5)$}" $ap^{2N_1+N_2+2N_3}$
\item"{$(1,2,6)$}" $-a^2\x p^{2N_3}$
\item"{$(1,3,4)$}" $-bp^{2N_1+N_2+2N_3} +\ab a^2p^{N_1+N_2+2N_3}$
\item"{$(1,3,5)$}" $p^{3N_1+N_2+2N_3}+a^2\x p^{2N_3}-a^2\aaa 
p^{N_1+N_2+3N_3}$
\item"{$(1,3,6)$}" $a\x p^{N_1+2N_3} +a^2(\ac\pnb -c)p^{N_1+2N_3}$
\item"{$(1,4,5)$}" $\x bp^{N_1+N_2+N_3}+\ab p^{3N_1+2N_2+N_3}-\aaa b 
p^{2N_1+2N_2+N_3} $
\item"{$(1,4,6)$}" $(\ac\pnb -c)bp^{2N_1+N_2+N_3}+\ab a \x p^{N_1+N_2+N_3}$
\item"{$(1,5,6)$}" $(\ac\pnb -c)p^{3N_1+N_2+N_3}+a\x^2p^{N_3}-a\x 
p^{N_1+N_2+N_3}$
\item"{$(2,3,4)$}" $-a^2\x p^{2N_3} $
\item"{$(2,3,5)$}" $-a^2p^{N_1+N_2+2N_3}$
\item"{$(2,3,6)$}" $0$
\item"{$(2,4,5)$}" $-acp^{2N_1+N_2+N_3}$
\item"{$(2,4,6)$}" $-a\x^2\pnc$
\item"{$(2,5,6)$}" $a\x p^{N_1+N_2+N_3}$
\item"{$(3,4,5)$}" $cp^{3N_1+N_2+N_3}+a\x^2 \pnc -\aaa a\x 
p^{N_1+N_2+N_3}-a\ab p^{2(N_1+N_2)+N_3}$
\item"{$(3,4,6)$}" $\x ^2\pnac +a\x (\ac\pnb-c)\pnac$
\item"{$(3,5,6)$}" $(\ac\pnb-c)ap^{2N_1+N_2+N_3}$
\item"{$(4,5,6)$}" $-\x^3+\aaa \x^2\pnab +\ab \x p^{2(N_1+N_2)}-c^2 
p^{2N_1}\pnab +\ac c p^{2N_1}\pnab$
\endroster

So for example, to get rid of the non-monomial expression
$$
(\ac\pnb -c)p^{3N_1+N_2+N_3}+a\x^2p^{N_3}-a\x p^{N_1+N_2+N_3}
$$
coming from the determinant $(1,5,6)$ we apply the condition
$$
\min \{v(X+Y),v(X)\}=\min \{v(X),v(Y)\}
$$
to the determinants $(2,4,6)$ and $(2,5,6)$ to eliminate from consideration 
the terms $a\x^2p^{N_3}$ and $a\x\pnabc.$ The same process allows us to 
neglect all non-monomial expressions except that arising from $(4,5,6).$

\vskip 0.5in
\heading{Bibliography}
\endheading

\
\newline
[Be] A Beauville, {\it Determinantal hypersurfaces,} Michigan Math. J. {\bf
48} (2000), 39--64.
\newline
[BM] Brumer, McGuinness, {\it The behaviour of the Mordell-Weil Group of 
elliptic curves,} Bull. of the AMS {\bf 23} (1990), 375--382.
\newline
[D] J Denef, {\it The rationality of the Poincare series associated to the
$p$-adic points on a variety, } Invent. Math. {\bf  77} (1984), 1--23.
\newline
[DM] J Denef, D Meuser, {\it A functional equation of Igusa\rq s local zeta 
function,} Amer. J. Math. {\bf 113} (1991), 1135-1152.
\newline
[Di] L Dickson, {\it Determination of all general homogeneous polynomials
expressible as determinants  with linear elements,} Trans. Amer. Math. Soc.
{\bf 22} (1921), 167--179.
\newline
[duS1] M du Sautoy, {\it Finitely generated groups, $p$-adic analytic groups
and
Poin- \linebreak care series,} Ann. Math. {\bf 137} (1993), 639--670.
\newline
[duS2] M du Sautoy, {\it The zeta function of $\frak {sl}_2(\Bbb Z)$,}
Forum Mathematicum {\bf 12} (2000), 197--221.
\newline
[duS3] M du Sautoy, {\it Counting subgroups in nilpotent groups and points
in elliptic curves, }  MPI preprint 2000--86.
\newline
[duS4] M du Sautoy, {\it A nilpotent group and its elliptic curve:
non-uniformity of local zeta functions of groups, } MPI preprint 2000--85.
To appear in Israel J. Math. {\bf 126}.
\newline
[duS5] M du Sautoy, {\it Counting $p$-groups and subgroups of nilpotent
groups,} Extrait des Publ. Math. {\bf 92} (2000) 63--112.
\newline
[duSG1] M du Sautoy, F J Grunewald, {\it Analytic properties of zeta
functions and subgroup growth, } Ann. of Math.
{\bf 152} (2000), 793--833.
\newline
[duSG2] M du Sautoy, F J Grunewald, {\it Uniformity for 2-generator free
nilpotent groups,} in preparation.
\newline
[duSL] M du Sautoy, F Loeser, {\it Motivic zeta functions for infinite
dimensional Lie algebras, } Ecole Polytechnique preprint series 2000-12.
\newline
[GSS] F J Grunewald, D Segal, G C Smith, {\it Subgroups of finite index in
nilpotent groups, } Invent. Math. {\bf 93} (1988), 185--223.
\newline
[I] J Igusa, {\it An introduction to the theory of local zeta functions, }
Studies in Adv. Math., AMS, {\bf 14}
2001.
\newline
[M1] J S Milne, {\it Lecture notes on elliptic curves,} Univ. of  Michigan.
\newline
[M2] J S Milne, {\it Lecture notes on abelian varieties, } Univ. of
Michigan.
\newline
[Se] D Segal, {\it Polycyclic groups,} Cambridge Univ. Press, Cambridge,
1983.
\newline
[Sil] J H Silverman, {\it The arithmetic of elliptic curves,} GTM {\bf 106},
Springer-Verlag Publ., 1986.
\newline
[W] T Womack, computer search, Nottingham, 2002.
\newline
[V] C Voll, PhD Thesis, Cambridge, 2002.

\enddocument